\documentclass[12pt]{article}
\usepackage{amsmath,amssymb,amscd,color}
\newcommand{\bq}{\begin{equation}}
\newcommand{\eq}{\end{equation}}

\newcommand{\R}{{ \mathbb{R}  }}

\begin{document}
\bibliographystyle{plain}

\newtheorem{defn}{Definition}
\newtheorem{lemma}[defn]{Lemma}
\newtheorem{proposition}{Proposition}
\newtheorem{theorem}[defn]{Theorem}
\newtheorem{cor}{Corollary}
\newtheorem{remark}{Remark}
\numberwithin{equation}{section}

\def\Xint#1{\mathchoice
   {\XXint\displaystyle\textstyle{#1}}%
   {\XXint\textstyle\scriptstyle{#1}}%
   {\XXint\scriptstyle\scriptscriptstyle{#1}}%
   {\XXint\scriptscriptstyle\scriptscriptstyle{#1}}%
   \!\int}
\def\XXint#1#2#3{{\setbox0=\hbox{$#1{#2#3}{\int}$}
     \vcenter{\hbox{$#2#3$}}\kern-.5\wd0}}
\def\ddashint{\Xint=}
\def\dashint{\Xint-}
\def\aint{\Xint\diagup}

\newenvironment{proof}{{\bf Proof.}}{\hfill\fbox{}\par\vspace{.2cm}}
\newenvironment{pfthm1}{{\par\noindent\bf
            Proof of Theorem \ref{thm1}. }}{\hfill\fbox{}\par\vspace{.2cm}}
\newenvironment{pfthm0}{{\par\noindent\bf
            Proof of Theorem \ref{global}. }}{\hfill\fbox{}\par\vspace{.2cm}}
\newenvironment{pfthm3}{{\par\noindent\bf
Proof of Theorem 3. }}{\hfill\fbox{}\par\vspace{.2cm}}
\newenvironment{pfthm4}{{\par\noindent\bf
Proof of Theorem \ref{theorem3}. }}{\hfill\fbox{}\par\vspace{.2cm}}
\newenvironment{pfthm5}{{\par\noindent\bf
Proof of Theorem 5. }}{\hfill\fbox{}\par\vspace{.2cm}}
\newenvironment{pflemsregular}{{\par\noindent\bf
            Proof of Lemma \ref{sregular}. }}{\hfill\fbox{}\par\vspace{.2cm}}

\title{ On the Geometric  Regularity Conditions  for   the 3D Navier-Stokes Equations }
\author{$\mbox{Dongho Chae}^{(*)}\mbox{ and Jihoon Lee}^{(\dagger)}$\\
 Department of Mathematics;\\
 Chung-Ang University\\
 Seoul 156-756; Korea\\
 e-mail : ($*$)dchae@cau.ac.kr,\quad ($\dagger$)jhleepde@cau.ac.kr}

\date{}

\maketitle
\begin{abstract}
We prove geometrically improved version of Prodi-Serrin type blow-up criterion.  Let $v$ and $\omega$ be the velocity and the vorticity of solutions to the 3D Navier-Stokes equations and  denote  $\{f\}_+=\max\{f, 0\}$ , $Q_T=\Bbb R^3\times (0, T)$.  If  $\left\{\left(  v \times \frac{\omega}{|\omega|}   \right)\cdot \frac{\Lambda^{\beta}v}{|\Lambda^{\beta}v|}\right\}_+ \in L^{\gamma, \alpha}_{x,t} (Q_T)$ with $3/\gamma +2/\alpha \leq 1$ for some $\gamma >3$ and $1 \leq \beta \leq 2$, then the local smooth solution $v$ of the Navier-Stokes equations on $(0,T)$ can be continued to $(0, T+\delta)$ for some $\delta >0$.  We also prove localized version of a special case of this. Let $v$ be a suitable weak solution to the Navier-tokes equations in a space-time domain containing $z_0= (x_0, t_0)$,  let $Q_{z_0, r}=B_{x_0, r} \times (t_0-r^2, t_0)$ be a parabolic cylinder in the domain. We show that if either $\left\{\left( v \times \frac{\omega}{|\omega|}\right) \cdot \frac{\nabla \times \omega}{|\nabla \times \omega|}\right\}_{+} \in L^{\gamma, \alpha}_{x,t}(Q_{z_0, r})$ with $\frac{3}{\gamma}+\frac{2}{\alpha} \leq 1$, or $\left\{\left(\frac{v}{|v|} \times \omega\right) \cdot \frac{\nabla \times \omega}{|\nabla \times \omega|}\right\}_{+} \in L^{\gamma, \alpha}_{x,t}(Q_{z_0, r})$ with $\frac{3}{\gamma}+\frac{2}{\alpha} \leq 2$, ($\gamma \geq 2$, $\alpha \geq 2$),  then $z_0$ is a regular point for $v$. This improves previous local regularity criteria for the suitable weak solutions.
\end{abstract}

\section{Introduction}
 \setcounter{equation}{0}

\indent
 In this brief note we consider three-dimensional incompressible Navier-Stokes equations in a domain $\Omega \subset {\mathbb{R}}^3$ :
\[
(NS)\left\{\begin{array}{ll}
 v_{t}+(v \cdot \nabla) v =-\nabla p +\Delta v, \quad& (x,t) \in \Omega \times (0,T)\\
\nabla \cdot v=0,\qquad& (x,t) \in \Omega \times (0,T)\\
v(x,0)=v_0(x),\qquad & x \in \Omega
\end{array}\right.\]
 where $v=(v_1, v_2, v_3) $ is the flow velocity and $p$ is the scalar pressure, respectively. The initial data $v_0$ satisfies
\[ \nabla \cdot v_0 =0.
\]
It is well known that the first equations of  $(NS)$ can be rewritten as following equivalent form:
\begin{equation}\label{13}
v_t -v \times \omega =-\nabla \left( p+\frac{|v|^2}{2}   \right)+\Delta v,
\end{equation}
where $\omega=\nabla \times v$ is the vorticity vector field. The global in time existence of a smooth solution to the system (NS) is one of the outstanding open problems in mathematical fluid mechanics.
On the other hand,  the global in time existence of weak solution(Leray-Hopf weak solution) was proved first by Leray\cite{leray}.
There are numerous conditional regularity results of weak solutions by imposing the integrability conditions on the velocity or vorticity using scaling invariant function space for weak solutions to (NS) (see \cite{Bei1, Bei2, CKL, ESS, Giga, lady, ohyama, prodi, serrin} and references therein).
Besides the so-called Prodi-Serrin type regularity conditions, there are many studies on the geometric regularity conditions by imposing alignment of the direction of the vorticity (see \cite{Bei-Ber, Ber, Ber-Cor, Chae2, Const-Feff, Gruj1, Gruj2} and references therein).
Among the previous results, Chae\cite{Chae1} obtained local regularity criterion by imposing scaling invariant integrability conditions on $ v \times \frac{\omega}{|\omega|}$ or $\omega \times \frac{v}{|v|}$ which is a refinement of other Prodi-Serrin type condition on $v$ and $\omega$. 
On the other hand, Lee\cite{lee} obtained regularity by assuming the smallness of the volume of the parallelepiped which is defined by the unit vectors $\frac{v}{|v|}$, $\frac{\omega}{|\omega|}$ and $\frac{\nabla \times \omega}{|\nabla \times \omega|}$.\\
 We define nonlocal operator $\Lambda =(-\Delta)^{\frac12}$ as $\Lambda^{\beta} f= (-\Delta)^{\frac{\beta}{2}} f ={\mathcal{F}}^{-1}(|\xi|^{\beta} {\mathcal{F}}{f}(\xi))$ where ${\mathcal{F}}$ denotes a Fourier transform on $\R^3$. We use a mixed type norms for $Q_{T}= \R^3 \times (0, T)$ :
\[
\| v\|_{ L^{\gamma, \alpha}_{x,t}(Q_{T})}:= \left\|  \| v(\cdot, t)\|_{L^{\gamma}_{x}(\R^3) }  \right\|_{L^{\alpha}_{t}(0, T)},\quad 1\leq \alpha, \gamma \leq \infty.
\]
We also denote $\{f\}_+(x):=\max \{ f(x),\, 0\}$. Also direction fields $\frac{\omega}{|\omega|}$ and $\frac{\Lambda^{\beta} v}{|\Lambda^{\beta}v|}$ are set to be zero when $\omega(x,t)=0$ and $\Lambda^{\beta}v(x,t)=0$, respectively.\\ First, we consider Prodi-Serrin type blow-up criterion in terms of some triple product, which improves the previous criterion of \cite{serrin}. We consider only $\Omega=\R^3$ case for  simplicity.
\begin{theorem}\label{global}
Let $v$ be a local in time reguar solution of the Navier-Stokes equations (NS) in $Q_{T}:=\R^3 \times (0, T)$ with $v_0 \in H^{\frac{1}{2}}(\R^3)$.   Then, we have,
\begin{description}
\item[(i)] if $v$, $\omega:=\nabla \times v$ and $\Lambda^{\beta} v$ satisfies that, for an absolute constant $\epsilon_0$ and some $\beta \in [1,\,2]$, 
\begin{equation}\label{global-small-cond}
\left\|\left\{  \left( v \times \frac{\omega}{|\omega|}   \right)\cdot \frac{\Lambda^{\beta} v}{|\Lambda^{\beta} v|}\right\}_{+}\right\|_{L^{3, \infty}_{x,t}(Q_{T})}\leq \epsilon_0,
\end{equation}
then  a regular solution $v$ exists beyond $T$, that is, $v \in C([0,\, T+\delta) ; H^{\frac12}(\R^3))$ for some $\delta >0$.
\item[(ii)] $v$ blows up at $T_{*}$, which is a finite maximal time of local in time  smooth solution to (NS), namely,
\[
\limsup_{ t \nearrow T_{*}} \| v(t)  \|_{H^{m}} =\infty,\qquad \forall m \geq \frac12,
\]
if and only if for all 
 $\gamma \in (3, \infty]$ and $\alpha \in [2, \infty]$ with $3/\gamma +2/\alpha \leq 1$ and all $\beta \in [1.\, 2]$
\begin{equation}\label{global-v-cond}
 \left\|\left\{\left( v \times \frac{\omega}{|\omega|}   \right)\cdot \frac{\Lambda^{\beta} v}{|\Lambda^{\beta} v|}\right\}_{+} \right\|_{ L^{\gamma, \alpha}_{x,t}(Q_{T})}=\infty.
\end{equation}
\end{description}
\end{theorem}
\begin{remark}
From the standard local in time existence theory of Navier-Stokes equations, $v(t) \in H^{m}(\R^3)$ for any $m \in {\mathbb{N}}$ and  $t \in (0,\, T_{*})$ where $T_{*}$ is a possible blow up time of local $H^{\frac12}$-solution. Therefore, any derivatives in Theorem \ref{global} are well-defined pointwise and $\Lambda^{\beta} v$ can be used as a test function.
\end{remark}
Since $\Lambda$ is a nonlocal operator, it does not seem easy to obtain local regularity criterion for Theorem \ref{global}. But for the case $\beta=2$, we can obtain a local regularity criterion for the triple product including $v$, $\omega$ and $-\Delta v = \nabla \times \omega$. \\
Our goal in this paper is to prove local regularity criterion by imposing integrability conditions on the triple product $\left( v \times \frac{\omega}{|\omega|}\right)\cdot \frac{\nabla \times \omega}{|\nabla \times \omega|}$ or $\left(\frac{v}{|v|} \times    \omega \right)\cdot \frac{\nabla \times \omega}{|\nabla \times \omega|}$. This improves both of the results in \cite{Chae1, lee} as well as the previous Prodi-Serrin type conditions.
For the local analysis of a weak solution the notion of  suitable weak solution is useful as shown in  the partial regularity results (see \cite{CKN} and \cite{scheffer}).
Let $Q_{T}:= \Omega \times (0, T),$ For a point $z=(x,t)\in Q_{T}$, we denote
\[
B_{x,r} := \{ y \in {\mathbb{R}}^3 \, :\, |y-x|<r\},\quad Q_{z,r}=B_{x,r} \times (t-r^2, t).
\]
We also use the mixed space-time norms :
\[
\| v\|_{ L^{\gamma, \alpha}_{x,t}(Q_{z,r})}:= \left\|  \| v(\cdot, t)\|_{L^{\gamma}_{x}(B_{x,r}) }  \right\|_{L^{\alpha}_{t}(t-r^2, t)},\quad 1\leq \alpha, \gamma \leq \infty.
\]
We state the definition of a suitable weak solution to (NS) for local analysis.
\begin{defn}
A pair $(v,p)$ of measurable functions is a suitable weak solution of (NS) if the following conditions are satisfied :
\begin{description}
\item[(i)] $ v \in L^{\infty}(0, T; L^2(\Omega))\cap L^2(0, T; W^{1,2}(\Omega))$, $p \in L^{\frac32}(Q_T))$.
\item[(ii)] The pair $(v,p)$ satisfies (NS) in the sense of distribution.
\item[(iii)] The pair $(v,p)$ satisfies the local energy inequality,
\[
\int_{\Omega} |v(x,t)|^2 \phi(x,t) dx +2\int_0^t \int_{\Omega} |\nabla v(x, \tau)|^2 \phi (x, \tau) dx d\tau
\]
\[
\leq \int_0^t \int_{\Omega} \left(  |v|^2 (\partial_t \phi +\Delta \phi)+(|v|^2 +2p) v \cdot \nabla \phi   \right)dx d\tau
\]
for almost all $t \in (0, T)$ and all nonnegative scalar test function $\phi \in C_0^{\infty}(Q_{T})$.
\end{description}
\end{defn}
We say that a weak solution is regular at $z$, if $v$ is bounded in $Q_{z,r}$ for some $r>0$. This point $z$ is called a regular point. 

Below we use extended definitions of the directional fields $v(x,t)/|v(x,t)|$, $\omega(x,t)/|\omega(x,t)|$ and $\nabla \times \omega(x,t)/|\nabla \times \omega(x,t)|$, which are set to zero whenever $v(x,t)=0$, $\omega(x,t)=0$ and $ \nabla \times \omega(x,t)=0$, respectively. 
\begin{theorem}\label{thm1}
Let $z_0 =(x_0, t_0) \in Q_{T}$ with $\bar{Q}_{z_0,r}  \subset Q_{T}$, and $(v,p)$ be a suitable weak solution of (NS) in $Q_{T}$ with the vorticity $ \omega =\nabla \times v$, where the derivatives are in the sense of distribution. Suppose $v$ and $\omega$ satisfy one of the following conditions :
\begin{description}
\item[(i)] There exists an absolute constant $\epsilon_0$ such that
\begin{equation}\label{CL-smallness-cond}
\left\|\left\{  \left( v \times \frac{\omega}{|\omega|}   \right)\cdot \frac{\nabla \times \omega}{|\nabla \times \omega|}\right\}_{+}\right\|_{L^{3, \infty}_{x,t}(Q_{z_0,r})}\leq \epsilon_0.
\end{equation}
\item[(ii)] There exists $\gamma \in (3, \infty]$ and $\alpha \in [2, \infty]$ with $3/\gamma +2/\alpha \leq 1$ such that
\begin{equation}\label{CL-v-cond}
 \left\{\left( v \times \frac{\omega}{|\omega|}   \right)\cdot \frac{\nabla \times \omega}{|\nabla \times \omega|}\right\}_{+} \in L^{\gamma, \alpha}_{x,t}(Q_{z_0,r}).
\end{equation}
\item[(iii)] There exists $\gamma \in [2, \infty]$ and $\alpha \in [2, \infty]$ with $3/\gamma+2/\alpha \leq 2$ such that
\begin{equation}\label{CL-omega-cond}
 \left\{\left( \frac{v}{|v|} \times  \omega  \right)\cdot \frac{\nabla \times \omega}{|\nabla \times \omega|} \right\}_{+}\in L^{\gamma, \alpha}_{x,t}(Q_{z_0,r}).
\end{equation}
Then, $z_0$ is a regular point.
\end{description}
\end{theorem}
\begin{remark} We note that there are many physical flows, including Beltrami flows (see \cite{Const-Maj}), for which the triple product vanishes.  The above theorem says intuitively that  even if the flow is far from the Beltrami flows,  
if the projection of the vector $\nabla\times \omega $ on the plane spanned by $v$ and $\omega$ is  ``controllable" in a local space-time region, then the flow is smooth in that region.
\end{remark}
\begin{remark}
In \cite{lee}  it was proved that if there exists an absolute constant $\epsilon_0$ such that
\begin{equation}
\left\|  \left( \frac{v}{|v|} \times \frac{\omega}{|\omega|}   \right)\cdot \frac{\nabla \times \omega}{|\nabla \times \omega|}\right\|_{L^{\infty, \infty}_{x,t}(Q_{z_0,r})}\leq \epsilon_0,
\end{equation}
then $z_0$ is a regular point. As an easy consequence of (i) of Theorem \ref{thm1}, we can have, for $b>0$, that if there exists an absolute constant $\epsilon_0$ such that
\begin{equation}
\left\|  |v|^{b}\left\{\left( \frac{v}{|v|} \times \frac{\omega}{|\omega|}   \right)\cdot \frac{\nabla \times \omega}{|\nabla \times \omega|}\right\}_{+}\right\|_{L^{\frac{3}{b}, \infty}_{x,t}(Q_{z_0,r})}\leq \epsilon_0,
\end{equation}
then $z_0$ is a regular point. Hence, the result in \cite{lee} is a special case of Theorem \ref{thm1} as $b \rightarrow 0+$.
\end{remark}
\begin{remark}
Theorem \ref{thm1} (i) and (ii) can be considered as  improvements of Theorem 1.1 (i) and (ii) in \cite{Chae1}. But Theorem \ref{thm1} (iii) can extend Theorem 1.1 (iii) of \cite{Chae1} only on the range $ \gamma \in [2,3]$ due to the technical difficulties. In order to extend Theorem 1.1 (iii) of \cite{Chae1} to the triple product on the range $\gamma \in (\frac32, 2) \cup (3, \infty]$, it seems necessary to develop different methods. 
\end{remark}
\section{Proof of the Main Theorems}
First, we prove Theorem \ref{global} by using standard a priori estimates.\\

\begin{pfthm0}
Let $T_{*}$ be a maximal time of local existence of $H^{\frac{1}{2}}$ solution. Multiplying $\Lambda^{\beta} v$ on the both sides of \eqref{13} and integrating over $\R^3$, we have, for $t<T_{*}$
\[
\frac12 \frac{d}{dt} \| \Lambda^{\frac{\beta}{2}} v \|_{L^2}^2+\| \nabla \Lambda^{\frac{\beta}{2}}v\|_{L^2}^2= \int_{\R^3} ( v\times \omega)\cdot \Lambda^{\beta} v dx
\]
\[
\leq \int_{\R^3} \left\{ \left(v \times \frac{\omega}{|\omega|} \right) \cdot \frac{\Lambda^{\beta} v}{|\Lambda^{\beta} v|} \right\}_{+}|\omega|\, |\Lambda^{\beta} v|\, dx
\]
\[
\leq \left\|\left\{ \left(v \times \frac{\omega}{|\omega|} \right) \cdot \frac{\Lambda^{\beta} v}{|\Lambda^{\beta} v|} \right\}_{+}\right\|_{L^{\gamma}} \| \omega \|_{L^p}\| \Lambda^{\beta} v\|_{L^q}:=I,
\]
where $p$ and $q$ safisfies $\frac{1}{p}+\frac{1}{q} =\frac{\gamma-1}{\gamma}$, $p \in [\frac{6}{5-\beta},\, \frac{6}{3-\beta}]$ and $q\in [\frac{6}{3+\beta}, \,\frac{6}{1+\beta}]$.\\
By the interpolation inequality, we have
\[
\| \omega \|_{L^p} \leq C \| \Lambda^{\frac{\beta}{2}} v \|_{L^2}^{\frac{3}{p}+\frac{\beta}{2}-\frac{3}{2}} \|\nabla \Lambda^{\frac{\beta}{2}}v\|_{L^2}^{\frac52-\frac{3}{p}-\frac{\beta}{2}} 
\]
and
\[
\| \Lambda^{\beta} v \|_{L^q} \leq C\| \Lambda^{\frac{\beta}{2}} v \|_{L^2}^{\frac{3}{q}-\frac{\beta}{2}-\frac{1}{2}} \|\nabla \Lambda^{\frac{\beta}{2}}v\|_{L^2}^{\frac32-\frac{3}{q}+\frac{\beta}{2}} .
\]
Then we can estimate $I$ as
\[
I \leq C\left\|\left\{ \left(v \times \frac{\omega}{|\omega|} \right) \cdot \frac{\Lambda^{\beta} v}{|\Lambda^{\beta} v|} \right\}_{+}\right\|_{L^{\gamma}}\| \Lambda^{\frac{\beta}{2}} v\|_{L^2}^{\frac{\gamma-3}{\gamma}}\| \nabla \Lambda^{\frac{\beta}{2}} v\|_{L^2}^{\frac{\gamma+3}{\gamma}}.
\]
We first assume the condition (i) of Theorem \ref{global} holds true. Then we have
\[
\frac12 \frac{d}{dt} \| \Lambda^{\frac{\beta}{2}} v \|_{L^2}^2+\left[ 1-C \left\|\left\{ \left(v \times \frac{\omega}{|\omega|} \right) \cdot \frac{\Lambda^{\beta} v}{|\Lambda^{\beta} v|} \right\}_{+}\right\|_{L^{3}} \right]\| \nabla \Lambda^{\frac{\beta}{2}}v\|_{L^2}^2\leq 0
\]
If $\epsilon_0 < \frac{1}{C}$, then $ v \in L^{\infty} (0, T_{*}; H^{\frac{\beta}{2}}(\R^3))$. By the standard continuation argument, we have $v \in C((0, T_{*} +\delta); H^{\frac{1}{2}}(\R^3))$ for some $\delta>0$.\\
Next, we assume the condition (ii) of Theorem \ref{global} holds true. By Young's inequality, we have
\[
I \leq  C\left\|\left\{ \left(v \times \frac{\omega}{|\omega|} \right) \cdot \frac{\Lambda^{\beta} v}{|\Lambda^{\beta} v|} \right\}_{+}\right\|_{L^{\gamma}}^{\frac{2\gamma}{\gamma-3}}\| \Lambda^{\frac{\beta}{2}} v\|_{L^2}^2 +\frac12 \|  \nabla \Lambda^{\frac{\beta}{2}} v\|_{L^2}^2.
\]
Therefore, we obtain
\[
\frac{d}{dt} \| \Lambda^{\frac{\beta}{2}} v \|_{L^2}^2+\| \nabla \Lambda^{\frac{\beta}{2}}v\|_{L^2}^2 \leq C\left\|\left\{ \left(v \times \frac{\omega}{|\omega|} \right) \cdot \frac{\Lambda^{\beta} v}{|\Lambda^{\beta} v|} \right\}_{+}\right\|_{L^{\gamma}}^{\frac{2\gamma}{\gamma-3}}\| \Lambda^{\frac{\beta}{2}} v\|_{L^2}^2.
\]
By Gronwall's inequality, we have
\[
\sup_{t \in [0, T_{*})} \| \Lambda^{\frac{\beta}{2}} v (t)\|_{L^2}^2\leq \|v_0 \|_{H^{\frac{\beta}{2}}}^2 \exp \left[ C\left\|\left\{ \left(v \times \frac{\omega}{|\omega|} \right) \cdot \frac{\Lambda^{\beta} v}{|\Lambda^{\beta} v|} \right\}_{+}\right\|_{L^{\gamma, \frac{2\gamma}{\gamma-3}}_{x,t}(Q_{T_{*}})}^{ \frac{2\gamma}{\gamma-3}}   \right].
\]
 Note that $\left\|\left\{ \left(v \times \frac{\omega}{|\omega|} \right) \cdot \frac{\Lambda^{\beta} v}{|\Lambda^{\beta} v|} \right\}_{+}\right\|_{L^{\gamma, \frac{2\gamma}{\gamma-3}}_{x,t}(Q_{T_{*}})}^{ \frac{2\gamma}{\gamma-3}} < \infty$ due to \eqref{global-v-cond}. Hence, $v \in C((0, T_{*} +\delta); H^{\frac{1}{2}}(\R^3))$ for some $\delta>0$. This concludes the proof.
\end{pfthm0}
Before proceding our proof, we recall the notion of an epoch of possible irregularity of the suitable weak solution of the Navier-Stokes equations. It is well known that for weak solution there exists a closed set $E \subset I=[0, T]$ such that solutions are regular on $I \setminus E$ and $1/2$-dimensional Hausdorff measure of $E$ is zero. Moreover, $E$ can be written as $\displaystyle{I\setminus \left\{\bigcup_{i \in {\mathcal{I}}} (\alpha_i, \beta_i)\right\}}$ where ${\mathcal{I}}$ is at most countable and $(\alpha_i, \beta_i)$ are disjoint open intervals in $[0, T]$. As in \cite{Galdi}, we call $\beta_i$ as an epoch of possible irregularity. We recall the following Lemma proved by Neustupa and Penel\cite{neu-pen} on the epoch of possible irregularity for suitable weak solutions.
\begin{lemma}\label{lem1}
Let $z_0=(x_0, t_0) \in Q_{T}$. Suppose $v$ is a suitable weak solution of the Navier-Stokes equations in $Q_{T}$ and $t_0$ be an epoch of possible irregularity. Then there exist positive numbers $\tau$, $r_1$ and $r_2$ with $r_1 < r_2$ such that the follwings are satisfied :
\begin{description}
\item[(a)] $\tau$ is sufficiently small so that $t_0$ is only one epoch of possible irregularity in time interval $[t_0-\tau, t_0]$.
\item[(b)] The closure $B_{x_0, r_2} \times (t_0-\tau, t_0)$ is contained in $Q_{T}$, i.e., $\bar{B}-{x_0, r-2} \times [t_0-\tau, t_0] \subset Q_{T}$.
\item[(c)] $((\bar{B}_{x_0, r_2} \setminus B_{x_0, r_1}) \times [t_0-\tau, t_0])\cap {\mathcal{S}}=\phi$, where ${\mathcal{S}}$ is the set of possible singular points of $v$.
\item[(d)] $v$, $v_t$, and $p$ are, together with all their space derivatives, continuous on $(\bar{B}_{x_0, r_2} \setminus B_{x_0, r_1}) \times [t_0-\tau, t_0].$
\end{description}
\end{lemma}

\begin{pfthm1}
First, we assume that $t_0$ is an epoch of possible irregularity for $v$ in $Q_{z_0, r}$. Suppose that $0<r_1<r_2<r$ and $r^2 <\tau$ are the positive numbers in Lemma \ref{lem1}. For simplicity, we denote $B_1 =B_{x_0, r_1}$ and $B_{2}=B_{x_0, r_2}$. We choose cut-off function $\varphi \in C_0^{\infty}(B_2)$ such that $\varphi=1$ on $B_1$ and set $u=\varphi v -V$ where $V \in C_0^2 (B_2 \setminus \bar{B}_1)$ satisfies $\mbox{div }V= ( v\cdot \nabla)\varphi$. We note that $(v \cdot \nabla)\varphi$ satisfied the compatibility condition :
\[
\int_{B_2 \setminus \bar{B}_1} (v\cdot \nabla)\varphi dx = \int_{\partial B_2}\varphi v \cdot n_2 dS -\int_{\partial B_1}  v \cdot n_1 dS=0,
\]
where $n_i$ is a unit outward normal vector to the sphere $\partial B_i$.
 Using Bogovski$\rm{\breve{i}}$'s Theorem(see \cite{Bogo} or \cite[Theorem III.3.1]{Galdi}), we can prove that there exists at least one $V$ satisfying above properties.  Then, by a straightforward calculation,  $u$ satisfies
\begin{equation}\label{eq-main}
u_t -\varphi v \times \omega +\nabla \left( \varphi \left( p +\frac{|v|^2}{2}\right)\right) -\Delta u =h , \qquad \mbox{div }u=0,
\end{equation}
where $h$ satisfies 
\[
h=-\frac{\partial V}{\partial t} +\left(   p +\frac{|v|^2}{2} \right) \nabla \varphi -v \Delta \varphi -2 (\nabla \varphi \cdot \nabla)v +\Delta V.
\]
We note that $h(\cdot, t)$ is sufficiently smooth and supported in the region $(\bar{B}_2\setminus B_1)$. 
Multiplying $-\Delta u$ on the both sides of \eqref{eq-main} and integrating, we have
\[
\frac12\frac{d}{dt} \| \nabla u \|_{L^2(B_2)}^2 + \| \Delta u \|_{L^2(B_2)}^2\]
\[ = \int_{B_2} v \times (\varphi \omega ) \cdot (\nabla \times \nabla \times (\varphi  v)) dx -\int_{B_2} v \times (\varphi \omega) \cdot \nabla\times( \nabla \times V) dx -\int_{B_2} \Delta u \cdot h dx
\]
\[
\leq   \int_{B_2} v \times (\varphi \omega ) \cdot (\nabla \times \nabla \times (\varphi  v)) dx   + C\| v\|_{L^2}^2 \| \varphi \omega \|_{L^2}^2\]
\[ +C\| D^2 V \|_{L^{\infty}(B_2)}^2 +\frac{1}{8} \| \Delta u \|_{L^2(B_2)}^2 +C\| h \|_{L^2(B_2)}^2
\]
\[
\leq C \int_{B_2} |v|^2 \, |\varphi \omega|\, |\nabla^2 \varphi| dx +C \int_{B_2} |v| \, |\varphi \omega|\, |\nabla \varphi|\, |\nabla v | dx+ \int_{B_2} \varphi^2\left((v \times   \omega) \cdot ( \nabla \times \omega)\right)_{+} dx 
\]
\[
  + C\| v\|_{L^2}^2 \| \varphi \omega \|_{L^2}^2
 +C\| D^2 V \|_{L^{\infty}(B_2)}^2 +\frac{1}{8} \| \Delta u \|_{L^2(B_2)}^2 +C\| h \|_{L^2(B_2)}^2
\]
\[
 := I_1+I_2+I_3+ C\| v\|_{L^2}^2 \| \varphi \omega \|_{L^2}^2
 +C\| D^2 V \|_{L^{\infty}(B_2)}^2 +\frac{1}{8} \| \Delta u \|_{L^2(B_2)}^2 +C\| h \|_{L^2(B_2)}^2.
\]

$I_1$ and $I_2$ can be easily estimated as follows :
\[
I_1 \leq C \int_{B_2} |v|^2 \, |\nabla \times u -\nabla \varphi \times v +\nabla \times V|\, |\nabla^2 \varphi| dx
\]
\[
\leq C\|v \|_{L^3(B_2\setminus B_1)}^2 \| \nabla u \|_{L^2(B_2)}^{\frac12}\| \Delta u \|_{L^2(B_2)}^{\frac12}+ C (\| v \|_{L^3(B_2\setminus B_1)}^3+1)
\]
\[
\leq C\| v\|_{L^3(B_2\setminus B_1)}^2 \| \nabla u \|_{L^2(B_2)}^2 + \frac18 \| \Delta u \|_{L^2(B_2)}^2 +C (\| v \|_{L^3(B_2\setminus B_1)}^3+1),
\]
and 
\[
I_2 \leq C \int_{B_2} |v| \, |\nabla \times u -\nabla \varphi \times v +\nabla \times V|\, |\nabla \varphi|\, |\nabla v| dx
\]
\[
\leq C(\| v\|_{L^3(B_2\setminus B_1)}\| \nabla u\|_{L^6(B_2)}+\| v\|_{L^4(B_2\setminus B_1)}^2 +\| v\|_{L^2(B_2\setminus B_1)}) \| \nabla v \|_{L^2(B_2)}
\]
\[
\leq \frac18 \| \Delta u \|_{L^2(B_2)}^2 +C(\| v\|_{L^3(B_2\setminus B_1)}^2+1) \| \nabla v \|_{L^2(B_2)}^2.
\]
Here, we note that 
\[
\| v\|_{L^3(B_2\setminus B_1)} \leq C,
\]
for some constant $C$ and all $t \in [t_0-r_2^2, t_0]$ due to the choice of $r_1$ and $r_2$ in Lemma \ref{lem1}.\\
Let us set $\kappa := \left\{\left( v \times \frac{\omega}{|\omega|}   \right)\cdot \frac{\nabla \times \omega}{|\nabla \times \omega|}\right\}_{+}$. Then $I_3$ can be estimated as 
\[
I_3 \leq \int_{B_2} \kappa \, |\varphi \omega|\, |\varphi \nabla \times \omega|\, dx
\]
\[
\leq \int_{B_2} \kappa\, | \nabla \times u-\nabla \varphi \times v+\nabla \times V|\, |\Delta u -\nabla \varphi \times \omega-\Delta \varphi v +\Delta V|\, dx
\]
\[
\leq C\int_{B_2} \kappa\, |\nabla u |\, |\Delta u| dx +C\int_{B_2} \kappa\, |g_1 |\, |\Delta u| dx+C\int_{B_2} \kappa\, |\nabla u |\, |g_2| dx
\]
\[
:=I_3^1+I_3^2+I_3^3,
\]
where we set $g_1=\nabla \varphi \times v-\nabla \times V$ and $g_2 =\nabla \varphi \times \omega+\Delta \varphi v -\Delta V$. Since $g_1$ and $g_2$ are smooth functions supported on $(B_2 \setminus \bar{B}_1) \times (t_0-\tau, t_0]$, we estimate 
\[
I_3^2 ,\, I_3^3 \leq C \| v\|_{L^2} \| g_1\|_{L^{\infty}} \| \Delta u \|_{L^2}+ C \| v\|_{L^2} \| g_2 \|_{L^{\infty}} \| \nabla u \|_{L^2} \leq C \| v\|_{L^2}^2 + C\| \nabla u \|_{L^2}^2 +\frac18 \| \Delta u \|_{L^2}^2.
\]
We first assume the condition of Theorem \ref{thm1} holds true. In this case,  we estimate
\begin{equation}
I_3^1 \leq  C \| \kappa \|_{L^3(B_2)} \| \nabla u \|_{L^6 (B_2)} \| \Delta u \|_{L^2(B_2)} \leq C_1 \epsilon_0 \| \Delta u \|_{L^2(B_2)}^2
\end{equation}
Combining all the estimates $I_1$, $I_2$, $I_3^1$, $I_3^2$ and $I_3^3$, we have
\begin{eqnarray}
 \frac{d}{dt} \| \nabla u \|_{L^2(B_2)}^2 + \| \Delta u \|_{L^2(B_2)}^2 & \leq&  2C_1 \epsilon_0 \|\Delta u \|_{L^2(B_2)}^2\nonumber\\
&&+ C(\| \nabla u \|_{L^2(B_2)}^2+\| \nabla v \|_{L^2(B_2)}^2+ \| h \|_{L^2}^2+1) \label{CL-21}
\end{eqnarray}
for $t \in (t_0-r_2^2, t_0]$, and an absolute constant $C_1$. If $C_1 \epsilon_0 <\frac12$, then integrating \eqref{CL-21} in time over $[t_0-r_2^2, t_0]$, we can obtain  $\nabla u \in L^{2, \infty}_{x,t}(Q_{z_0, r_2})$, and therefore $\nabla v \in L^{2, \infty}_{x,t}(Q_{z_0, r_1})$. Applying Corollary 2.1 in \cite{Chae1}, we conclude that $z_0$ is a regular point.\\
Next, we assume that the condition (ii) of Theorem \ref{thm1} holds true, and estimate
\begin{eqnarray}
I_3^1 &\leq& C \| \kappa \|_{L^{\gamma}(B_2)} \| \nabla u \|_{L^{\frac{2\gamma}{\gamma-2}}} \| \Delta u \|_{L^2} \nonumber\\
&\leq & C \| \kappa \|_{L^{\gamma}(B_2)} \| \nabla u \|_{L^2}^{\frac{\gamma-3}{\gamma}} \| \Delta u \|_{L^2}^{\frac{\gamma+3}{\gamma}} \nonumber \\
& \leq & C \| \kappa \|_{L^{\gamma}(B_2)}^{\frac{2\gamma}{\gamma-3}} \| \nabla u \|_{L^2}^2 +\frac18 \| \Delta u \|_{L^2}^2,
\end{eqnarray}
where we used the interpolation inequality,
\[
\| \nabla u \|_{L^{\frac{2\gamma}{\gamma-2}}} \leq C \| \nabla u \|_{L^2}^{1-\frac{3}{\gamma}} \| \Delta u \|_{L^2}^{\frac{3}{\gamma}},
\]
for $ \gamma \in (3, \infty]$. Since $\kappa \in  L^{\gamma, \alpha}_{x,t}(Q_{z_0,r_2})$ with $3/\gamma+2/\alpha \leq 1$ and $\gamma>3$, we have
\[
\| \kappa \|_{L^{\gamma, \frac{2\gamma}{\gamma-3}}_{x,t}(Q_{z_0, r_2})}^{\frac{2\gamma}{\gamma-3}} \leq  \| \kappa \|_{L^{\gamma, \alpha}_{x,t}(Q_{z_0, r_2})}^{\frac{2\gamma}{\gamma-3}} r_2^{\frac{2\gamma}{\gamma-3}(1-\frac{3}{\gamma}-\frac{2}{\alpha})}<\infty.
\]
Similarly to the previous case, we conclude that $z_0$ is a regular point for $v$ by Gronwall's inequality.\\
Let us set $\eta := \left\{\left(\frac{v}{|v|}  \times \omega   \right)\cdot \frac{\nabla \times \omega}{|\nabla \times \omega|}\right\}_{+}$. Then $I_3$ can be estimated as 
\[
I_3 \leq \int_{B_2} \eta \, |\varphi v|\, |\varphi \nabla \times \omega|\, dx
\]
\[
\leq \int_{B_2} \eta\, | u+ V|\, |\Delta u -\nabla \varphi \times \omega-\Delta \varphi v +\Delta V|\, dx
\]
\[
\leq C\int_{B_2} \eta\, | u |\, |\Delta u| dx +C\int_{B_2} \eta \, |V |\, |\Delta u| dx+C\int_{B_2} \eta \, | u |\, |g_2| dx
\]
\[
:=J_3^1+J_3^2+J_3^3,
\]
where we set $g_2 =\nabla \varphi \times \omega+\Delta \varphi v -\Delta V$. Since $V$ and $g_2$ are smooth functions supported on $(B_2 \setminus \bar{B}_1) \times (t_0-\tau, t_0]$, we estimate 
\[
J_3^2 + J_3^3 \leq C \| \nabla v\|_{L^2}^2 + C\| u \|_{L^2}^2+\frac18 \| \Delta u \|_{L^2(B_2)}^2.
\]
Now we assume (iii) of Theorem \ref{thm1} holds true, then we estimate $J_3^1$ as
\begin{eqnarray*}
J_3^1 &\leq& C \|\eta \|_{L^{\gamma}(B_2)} \| u \|_{L^{\frac{2\gamma}{\gamma-2}}} \| \Delta u \|_{L^2}
\\
&\leq& \left\{  \begin{array}{ll} C\|\eta \|_{L^{\gamma}(B_2)} \| u \|_{L^6}^{\frac{2\gamma-3}{\gamma} } \| \Delta u \|_{L^2}^{\frac{3}{\gamma}} \qquad &\mbox{if } 2 \leq \gamma \leq 3\\  C\|\eta \|_{L^{\gamma}(B_2)} \| u \|_{L^2}^{\frac{\gamma-3}{\gamma} }\| \nabla u \|_{L^2}^{\frac{3}{\gamma} } \| \Delta u \|_{L^2} \qquad &\mbox{if } \gamma >3   \end{array}\right.\\
&\leq& \left\{  \begin{array}{ll} C\|\eta \|_{L^{\gamma}(B_2)}^{\frac{2\gamma}{2\gamma-3}} \|\nabla  u \|_{L^2}^2+\frac18 \| \Delta u \|_{L^2}^{2} \qquad &\mbox{if } 2 \leq \gamma \leq 3\\  C\|\eta \|_{L^{\gamma}(B_2)}^2 \| \nabla u \|_{L^2}^{\frac{6}{\gamma} } +\frac18\| \Delta u \|_{L^2}^2 \qquad &\mbox{if }  \gamma  >3  \end{array}\right.
\end{eqnarray*}
Since $\eta \in  L^{\gamma, \alpha}_{x,t}(Q_{z_0,r_2})$ with $3/\gamma+2/\alpha \leq 2$, $\gamma\geq 2$ and $\alpha \geq 2$, we have
\[
\| \eta \|_{L^{\gamma, \frac{2\gamma}{2\gamma-3}}_{x,t}(Q_{z_0, r_2})}^{\frac{2\gamma}{2\gamma-3}} \leq  \| \eta \|_{L^{\gamma, \alpha}_{x,t}(Q_{z_0, r_2})}^{\frac{2\gamma}{2\gamma-3}} r_2^{\frac{2\gamma}{2\gamma-3}(2-\frac{3}{\gamma}-\frac{2}{\alpha})}<\infty,
\]
and
\[
\| \eta \|_{L^{\gamma, 2}_{x,t}(Q_{z_0, r_2})}^2 \leq \| \eta  \|_{L^{\gamma, \alpha}_{x,t}(Q_{z_0, r_2})}^2r_2^{\frac{2(\alpha-2)}{\alpha}}< \infty.
\]
Similarly to the previous case, we conclude that $z_0$ is a regular point for $v$ by Gronwall's inequality.\\
Next, we suppose that $t_0$ is a singular time which is not an epoch of possible irregularity. Then there exists a time $t^{*}$ in $(t_0-r^2, t_0)$ and $0<\tilde{r}_1<\tilde{r}_2<r$ such that $v$ is regular on $B_{x_0, \tilde{r}_2} \setminus B_{x_0, \tilde{r}_1} \times [t^{*}, t_0]$. Assume that $v$ is not regular on $B_{x_0, \tilde{r}_1} \times [t^{*}, t_0]$, then there exists $s \in (t^{*}, t_0]$ such that the suitable weak solution is regular on $B_{x_0, \tilde{r}_1} \times [t^{*},s)$ and singularity occurs at $(y,s) \in B_{x_0, \tilde{r}_1}\times \{s\}$. Then we take a local neighborhood of $(y,s)$ contained in $B_{x_0, \tilde{r}_2} \times [t^{*},s).$ Hence we can show $(y,s)$ is a regular point by the repetition of the above argument as in the case of the epoch of possible irregularity. It gives a contradiction to the assumption that $(y,s)$ is a singular point and hence $v$ is regular on $B_{x_0, \tilde{r}_1} \times [t^{*}, t_0]$. This completes the proof.
\end{pfthm1}

\section*{Acknowledgments}
This work was partially supported by NRF grants no. 2016R1A2B3011647.

\end{document}